\documentclass[a4paper,11pt,reqno]{amsart}
\usepackage{amsfonts}
\usepackage{amsmath}
\usepackage{amssymb}
\usepackage[dvips]{color}

\newtheorem{theorem}{Theorem}

\newtheorem{corollary}{Corollary}
\newtheorem{remark}{Remark}
\newtheorem{lemma}{Lemma}

\def\be{\begin{equation}}
\def\ee{\end{equation}}
\def\ben{\begin{displaymath}}
\def\een{\end{displaymath}}
\def\baa{\begin{eqnarray}}
\def\eaa{\end{eqnarray}}

\def\ba{\begin{array}}
\def\ea{\end{array}}
\def\la{\label}
\def\p{\partial}

\makeatletter
\@addtoreset{equation}{section}
\makeatother

\def\C{{\mathbb C}}
\def\R{{\mathbb R}}
\def\Z{{\mathbb Z}}
\def\Q{{\mathbb Q}}
\def\E{{\mathbb E}}
\def\Eb{\overline{{\mathbb E}}}
\def\P{{\mathbb P}}
\def\H{{\mathcal H}}
\def\Ht{\tilde{{\mathcal H}}}
\def\Hc{\check{{\mathcal H}}}
\def\Hb{\overline{{\mathcal H}}}
\def\M{{\mathcal M}}
\def\Mb{\overline{{\mathcal M}}}

\def\2x2{{\left(\!\!\begin{array}{cc}a&b\\c&d\\\end{array}\!\!\right)}}
\def\del{\partial}
\def\f{\frac}
\def\e{\epsilon} 
\def\D{\Delta}
\def\d{\delta}
\def\a{\alpha}
\def\deg{{\rm deg}}

\begin{document}
\title[Tau function and moduli of differentials]{Tau function and moduli of differentials}

\author{D.~Korotkin}
\address{Department of Mathematics and Statistics, Concordia University,
1455 de Maisonneuve West, Montreal, H3G 1M8  Quebec,  Canada}
\email{korotkin@mathstat.concordia.ca}
\author{P.~Zograf}
\address{Steklov Mathematical Institute, Fontanka 27, St. Petersburg 191023 Russia}
\email{zograf@pdmi.ras.ru}
\thanks{DK was partially supported by NSERC, FQRNT and CURC; 
PZ was partially supported by the RFBR grant 08-01-00379-a and 
by the President of Russian Federation grant NSh-2460.2008.1.}

\begin{abstract}
The tau function on the moduli space of generic holomorphic 1-differentials on complex algebraic 
curves is interpreted as a section of a line bundle on the projectivized Hodge bundle over the moduli space of stable curves. 
The asymptotics of the tau function near the boundary of the moduli space of generic 1-differentials is computed, and an
explicit expression for the pullback of the Hodge class on the projectivized Hodge bundle in terms of the tautological
class and the classes of boundary divisors is derived. This expression is used to clarify the geometric meaning of the 
Kontsevich-Zorich formula for the sum of the Lyapunov exponents associated with the Teichm\"uller flow on the Hodge bundle.
\end{abstract}

\maketitle


\section{Introduction}

Moduli spaces of holomorphic 1-differentials on complex algebraic curves arise 
in various areas of mathematics from algebraic geometry to completely integrable systems, holomorphic dynamics and ergodic theory.
Notably, they admit an ergodic $SL(2,\R)$ action that can be desribed as follows. Take the Hodge bundle on the moduli space
$\M_g$ of complex algebraic curves (the fibers of this bundle are the spaces of holomorphic 1-differentials on the corresponding
curve), and consider it as a real analytic space. The group $SL(2,\R)$ acts by linear transformations on the real and imaginary
parts of a holomorphic 1-form. The dynamics of this action has been extensively studied by many authors. It is closely related 
to billiards in rational polygons and to interval exchange maps, and its invariants admit a nice geometric interpretation (cf.
the pioneering work \cite{KZ} for more details).

The isomonodromic tau function on Hurwitz spaces has a straightforward analogue on moduli spaces of holomorphic differentials.  
It can be explicitly written in terms of the theta function and the prime form on the underlying curve and plays an
important role in the holomorphic factorization of determinats of flat Laplacians \cite{KK}. Here we follow
the approach of \cite{KKZ} to study the asymptotic behavior of this tau function and to compute its divisor. 
This allows us to express the pullback of the Hodge class on the projectivized Hodge bundle as a linear combination of 
the tautological class and the classes of boundary divisors. The obtained expression allows us to interpret geometrically 
 the Kontsevich-Zorich formula for the sum of the Lyapunov exponents of the 
$diag(e^t,e^{-t})$-action $(t\in\R)$ on the Hodge bundle over $\M_g$.

A few words about the structure of this paper. Section 2 contains some preliminaries on the moduli space of generic 1-differentials.
In Section 3 we define the tau function, give an explicit formula 
for it (Theorem 1), study its transformation properties and interpret it as a holomorphic section of a line bundle on the 
projectivized Hodge bundle over the moduli space $\M_g$.
Section 4 contains the main results of the paper: asymptotic formulae for the tau function near the boundary components (Theorem 2), 
and a formula for the pullback to the projectivized Hodge line bundle expressing it as a linear combination of the tautological class and the classes 
of boundary divisors (Theorem 3). In Section 5 we discuss the Kontsevich-Zorich formula for the sum of the Lyapunov exponents.

\section{Spaces of holomorphic 1-differentials}

Let $C$ be a smooth complex algebraic curve of genus $g$, and let
$\omega$ be a non-zero holomorphic 1-differential on $C$. We call a holomorphic 
differential {\em generic} if it has exactly $2g-2$ simple zeroes. 
Two pairs $(C_1,\omega_1)$ and $(C_2,\omega_2)$
are called {\em equivalent} if there exists
an isomorphism $h: C_1\rightarrow C_2$ such that $\omega_1=h^*(\omega_2)$.
The moduli space of generic pairs $(C,\omega)$ defined modulo this relation we denote by $\Ht_g^0$. 
Additionally we will consider an equivalence relation for
holomorphic 1-differentials on Torelli marked curves. A {\em Torelli marking} is a choice
of symplectic basis $\alpha=\{a_i, b_i\}_{i=1}^g$ in the first homology group
$H_1(C)$ of $C$. A curve $C$ together with a symplectic basis $\alpha$ will be
denoted by $C^\alpha$. We say that two pairs $(C_1^{\alpha_1},\omega_1)$ and $(C_2^{\alpha_2},\omega_2)$
are {\em Torelli equivalent} if there exists an isomorphism $h: C_1\rightarrow C_2$ such that
$\omega_1=h^*(\omega_2)$ and $h_*(\alpha_1)=\alpha_2$ elementwise.
The moduli space of pairs $(C^{\alpha},\omega)$ modulo the Torelli equivalence we denote by $\Hc_g^0$. 
The space $\Hc_g^0$ is a smooth non-compact complex manifold of dimension $4g-3$. The symplectic group
$Sp(2g,\Z)$ acts on $\Hc_g^0$ by changing Torelli marking, and $\Ht_g^0=\Hc_g^0/Sp(2g,\Z)$. Both $\Ht_g^0$
and $\Hc_g^0$ enjoy a natural action of $\C^*$ (by multiplication of $\omega$) that commutes with the 
action of $Sp(2g,\Z)$.

In the sequel we will also deal with holomorphic 1-differentials with degenerate zeroes. Let $\omega$ have
$r$ zeroes of multiplicities $m_1,\dots,m_r$ with $m_1+\dots+m_r=2g-2$. We call $\mu=(m_1-1,\dots,m_r-1)$ 
the {\em degeneracy type} of $\omega$ (we may omit all zero entries of $\mu$).
The moduli space of holomorphic 1-differentials of a fixed degeneracy type $\mu$ defined 
modulo the above equivalence (resp. Torelli equivalence) we denote by 
$\Ht_g^{\mu}$ (resp. by $\Hc_g^{\mu}$). According to \cite{KZ2}, these spaces 
are connected when $\omega$ has at least one simple zero (otherwise they may have up to 3 connected components).
Everything said in the previous paragraph about the action
of $Sp(2g,\Z)$ and $\C^*$ applies to the spaces $\Hc_g^{\mu}$ and $\Ht_g^{\mu}$ as well.
The dimension of these spaces is $4g-3-|\mu|$, where $|\mu|=\sum_{k=1}^r(\mu_k-1)$ is
the {\em total} degeneracy. 

We describe a natural completion of the space $\Ht_g^0$. Let $\M_g$ be the moduli space of smooth genus $g$ curves,
and let $\Mb_g$ be its Deligne-Mumford compactification. The boundary $\Mb_g-\M_g$ is the union of $[g/2]+1$
irreducible divisors $\D_0,\D_1,\dots,\D_{[g/2]}$, where $\D_0$ is the (closure of the) set of  irreducible
curves of arithmetic genus $g$ with one node, and $\D_j,\;j=1,\dots,[g/2],$ parametrizes reducible curves with 
components of genus $j$ and $g-j$.
Denote by $\E_g\to\M_g$ the Hodge bundle, where the fiber $\E_g|_C$ over a point represented by a curve $C$ is given 
by $\Omega_C^1$, the space of holomorphic 1-forms on $C$, up to the action of ${\rm Aut}(C)$. The Hodge bundle extends naturally 
to a bundle $\Eb_g\to\Mb_g$ (understood in the sense of orbifolds or algebraic stacks). The fiber
of $\Eb_g$ over a point represented by a reducible curve $C=C_1\cup C_2$ is given by $\Omega_{C_1}^1\oplus\Omega_{C_2}^1$,
whereas over an irreducible curve it is given by the vector space $\Omega_{C';p,q}^1$ 
of meromorphic 1-differentials on the normalization $C'\to C,\; g(C')=g-1,$ with at most simple poles at the preimages
$p,q$ of the node of $C$ with opposite residues. We have a sequence of inclusions $\Ht_g^0\hookrightarrow\E_g\hookrightarrow\Eb_g$, such that
the image of $\Ht_g^0$ is an open dense subset in $\Eb_g$. Note that we have a similar inclusion $\Ht_g^\mu\hookrightarrow\Eb_g$
for any degeneracy type $\mu$.

Denote by $\Hb_g=\P(\Eb_g)$ the projectivization of the Hodge bundle on $\Mb_g$. The space $\Hb_g$ is a smooth compact complex orbifold (smooth Deligne-Mumford stack)
of dimension $4g-4$, and the factor $\H_g=\Ht_g^0/\C^*$ is naturally included in $\Hb_g$ as an open dense subset. The complement
$\Hb_g-\H_g$ is the union of $[g/2]+2$ divisors:
\be
\Hb_g-\H_g=D_{\rm deg}\cup D_0\cup\dots\cup D_{[g/2]}\;.\label{compl}
\ee
Here the divisor
$D_{\rm deg}=\Hb_g^1$ is the closure in $\Hb_g$ of the locus $\Ht_g^1/\C^*$ of degenerate 1-differentials considered up to 
a constant factor, and $D_j=\pi^*(\D_j),\;j=0,\dots,[g/2],$ are the pullbacks of the boundary divisors $\D_j\subset\Mb_g$
via the natural projection $\pi:\Hb_g\to\Mb_g$.

Let $L\to\Hb_g$ be the tautological line bundle on $\Hb_g$ associated with the projection $(\Eb_g-\Mb_g)\to\P(\Eb_g)=\Hb_g$,
and put $\psi=c_1(L)\in {\rm Pic}(\Mb_g)\otimes\Q$. Denote by $\lambda=\pi^*(c_1(\Eb_g))$ the Hodge class in ${\rm Pic}(\Hb_g)\otimes\Q$,
that is, the pullback of the class $c_1(\Eb_g)\in {\rm Pic}(\Mb_g)\otimes\Q$ via the projection $\pi:\Hb_g\to\Mb_g$. We also put 
$\delta_{i}=[D_{i}]$ for $i\neq 1$ and $\delta_1=\f{1}{2}[D_1]$ in ${\rm Pic}(\Hb_g)\otimes\Q$.

\begin{lemma}
The rational Picard group ${\rm Pic}(\Hb_g)\otimes\Q$ of the space $\Hb_g$ is freely generated over $\Q$ by the classes
$\psi, \lambda, \delta_0,\dots,\delta_{[g/2]}$.
\end{lemma}
\begin{proof} By a result of \cite{AC}, the rational Picard group ${\rm Pic}(\Mb_g)\otimes\Q$ is freely generated 
by the classes
$\lambda_1, \D_0,\dots,\D_{[g/2]}$. We use the well-known fact that for a rank $n$ complex vector bundle $E\to M$ on 
a smooth complex variety $M$ one has
$$CH^*(\P (E))\cong CH^*(M)[\psi]/(\psi^n+c_1(E)\psi^{n-1}+\dots +c_n(E)),$$
where $\psi$ is the first Chern class of the tautological line bundle on $\P (E)$ (cf. \cite{Fu}, Example 8.3.4; here $CH^*$ stands
for the Chow ring). In particular, ${\rm Pic}(\P (E))\cong{\rm Pic}(M)\oplus \Z\psi$. The techniques of e.g. \cite{ELSV} allow to extend this statement (with rational coefficients) to the Hodge bundle $\Eb_g\to\Mb_g$. It then yields
${\rm Pic}(\Hb_g)\otimes\Q\cong({\rm Pic}(\Mb_g)\otimes\Q)\oplus\Q\psi$, i.e. ${\rm Pic}(\Hb_g)$ is freely generated 
by the classes $\psi, \lambda, \d_0,\dots,\d_{[g/2]}$.
\end{proof}

\section{Tau function}

The aim of this paper is to establish a non-trivial relation between the classes $\psi,\lambda,\d_0,\dots,\d_{[g/2]}$ and 
$\d_{\rm deg}=[D_\deg]$ in ${\rm Pic}(\Hb_g)\otimes\Q$ by explicitly computing the divisor of the  tau function of \cite{KK}.
 
For a Torelli marked curve $C^\alpha$, denote by $B(x,y)$ the 
{\em Bergman bidifferential}, that is, the unique symmetric meromorphic bidifferential 
on $C\times C$ with a quadratic pole of biresidue 1 on the diagonal and zero $a$-periods. 
Its $b$-periods
\be
\omega_i=\int_{b_i}B(\cdot\,, y)dy
\ee
are the {\em normalized holomorphic 1-differentials} on $C^\alpha$, that is,
\be
\int_{a_j}\omega_i=\delta_{ij},\quad\quad \int_{b_j}\omega_i=\Omega_{ij},\quad\quad i,j=1,\dots, g,
\ee
where the matrix $\Omega=\{\Omega_{ij}\}_{i,j=1}^g$ is the {\em period matrix} of $C^\alpha$.
In terms of local parameters $\zeta(x),\zeta(y)$ near the diagonal $\{x=y\}\in C\times C$, the bidifferential $B(x,y)$
has the expansion
\be
B(x,y)=\left(\f{1}{(\zeta(x)-\zeta(y))^2}+\f{S_B(\zeta(x))}{6}
+O((\zeta(x)-\zeta(y))^2)\right)d\zeta(x) d\zeta(y),
\la{asW}
\ee
where $S_B$ is a projective connection on $C$ called the {\em Bergman projective connection}.

Consider the non-linear differential operator  
$S_\omega=\frac{\omega''}{\omega}-\frac{3}{2}\left(\frac{\omega'}{\omega}\right)^2$ (that is, 
the Schwarzian derivative of the abelian integral $\int^x \omega$
with respect to a local parameter $\zeta$ on $C$). 
For a holomorphic 1-differential $\omega$, $S_\omega$ is a meromorphic projective connection on $C$, 
so that the difference $S_B-S_\omega$ is a meromorphic quadratic differential. 
Suppose that $\omega$ has $r$ zeroes $x_1,\dots,x_r$ of multiplicities $m_1,\dots,m_r$;
its degeneracy type is $\mu=(m_1-1,\dots,m_r-1)$ (this includes the case $\mu=0$).
Take the trivial line bundle on the space $\Hc_{g}^\mu$ and consider the connection
\be
d_B=d+\f{2}{\pi\sqrt{-1}}\,\sum_{i=1}^{2g+r-1}\left(\int_{s_i}\f{S_B-S_\omega}{\omega}\right)dz_i.
\la{Bercon}
\ee
Here $s_i=-b_i,\,s_{i+g}=a_i$ for $i=1,\dots,g$, and $s_{2g+k}$ is a small circle about $x_{k}$ for $k=1,\dots,r-1$,
whereas $z_i=\int_{a_i}\omega$, $z_{i+g}=\int_{b_i}\omega$ for $i=1,\dots,g$, and $z_{2g+k}=\int_{x_{2g-2}}^{x_{k}}\omega$ for $k=1,\dots,r-1$
($z_1,\dots,z_{2g+r-1}$ serve as local complex coordinates on $\Hc_g^\mu$, cf. \cite{KZ}). As it is shown
in \cite {KK}, this connection is flat. The tau function $\tau_\mu=\tau_\mu(C^\alpha, \omega)$ is locally defined as a horizontal
(covariant constant) section of the trivial line bundle on $\Hc_{g}^\mu$ with respect to $d_B$, that is,
\footnote{This tau function is the 24-th power of the Bergman tau function studied in \cite {KK}.} 
\be
d_B\log\tau_\mu=0.
\la{tau}
\ee

Let us now recall an explicit formula for the tau function $\tau_\mu$ derived in \cite {KK}.
Take a nonsingular odd theta characteristic $\delta$ and consider the corresponding theta function
$\theta[\delta](v; \Omega)$, where $v=(v_1,\dots,v_g)\in \C^g$. Put 
$$\omega_\delta=\sum_{i=1}^g\frac{\del\theta[\delta]}{\del v_i}\left(0; \Omega\right)\,\omega_i\,.$$
All zeroes of the holomorphic 1-differential $\omega_\delta$ have even multiplicities, and 
$\sqrt{\omega_\delta}$ is a well-defined holomorphic spinor on $C$. Following \cite{F}, consider 
the prime form \footnote{The prime form $E(x,y)$ is a canonical section of the line bundle on 
$C\times C$ associated with the diagonal divisor $\{x=y\}\subset C\times C$.}
\be
E(x,y)=\frac{\theta[\delta]\left(\int_x^y\omega_1,\dots,\int_x^y\omega_g; \Omega\right)}
{\sqrt{\omega_\delta}(x)\sqrt{\omega_\delta}(y)}.
\ee
To make the integrals uniquely defined, we fix $2g$ simple closed loops in the homology classes $a_i,b_i$
that cut $C$ into a connected domain, and pick the integration paths that do not intersect the cuts. 
The sign of the square root is chosen so that 
$E(x,y)=\frac{\zeta(y)-\zeta(x)}{\sqrt{d\zeta}(x)\sqrt{d\zeta}(y)}(1+O((\zeta(y)-\zeta(x))^2))$ as $y\rightarrow x$,
where $\zeta$ is a local parameter such that $d\zeta=\omega_\delta$.

We introduce local coordinates on $C$ that we call {\em natural} (or {\em distinguished}) with respect to $\omega$.
We take $\zeta(x)=\int_{x_1}^x\omega$ as a local coordinate on $C-\{x_1,\dots,x_r\}$, and choose $\zeta_k$  
near $x_k\in C$ in such a way that $\omega=d(\zeta_k^{m_k+1})=(m_k+1)\zeta_k^{m_k}d\zeta_k$, $k=1,\dots,r$. In terms of these coordinates we have 
$E(x,y)=\frac{E(\zeta(x),\zeta(y))}{\sqrt{d\zeta}(x)\sqrt{d\zeta}(y)}$, and we define
\begin{eqnarray*}
E(\zeta,x_k)&=&\lim_{y\rightarrow x_k}E(\zeta(x),\zeta(y))\sqrt{\frac{d\zeta_k}{d\zeta}}(y),\\
E(x_k,x_l)&=&\lim_{\stackrel{\scriptstyle x\rightarrow x_k}{y\rightarrow x_l}}E(\zeta(x),\zeta(y))
\sqrt{\frac{d\zeta_k}{d\zeta}}(x)\sqrt{\frac{d\zeta_l}{d\zeta}}(y)\,.
\end{eqnarray*}
Let ${\mathcal A}^x$ be the Abel map with the basepoint $x$, and let $K^x=(K^x_1,\dots,K^x_g)$ be the vector of Riemann constants
\be
K^x_i=\frac{1}{2}+\frac{1}{2}\Omega_{ii}-\sum_{j\neq i}\int_{a_i}\left(\omega_i(y)\int_x^y\omega_j\right)dy
\label{rc}\ee
(as above, we assume that the integration paths do not intersect the cuts on $C$).
Then we have ${\mathcal A}^x((\omega))+2K^x=\Omega Z+Z'$ for some $Z,Z'\in\Z^g$ . Now put
\be
\tau_\mu(C^\alpha, \omega)=
\f{\left(\left.\left(\sum_{i=1}^g\omega_i(\zeta)\frac{\partial}{\partial v_i}\right)^g
\theta(v;\Omega)\right|_{v=K^{\zeta}}\right)^{16}}{e^{4\pi\sqrt{-1}\langle\Omega Z+4K^{\zeta},Z\rangle}\;W(\zeta)^{16}}
\;\frac{\prod_{k<l}E(x_k,x_l)^{4m_k m_l}}
{\prod_k E(\zeta,x_k)^{8(g-1)m_k}}\;.
\label{taukk}
\ee
Here $\theta(v;\Omega)=\theta[0](v;\Omega)$ is the Riemann theta function, $v=(v_1,\dots,v_g)$,
and $W$ is the Wronskian of the normalized holomorphic differentials $\omega_1,\dots,\omega_g$ on 
$C^\alpha$. \footnote{The expression 
${\mathcal C}(\zeta)=\frac{1}{W(\zeta)}\left.\left(\sum_{i=1}^g\omega_i(\zeta)\frac{\partial}{\partial v_i}\right)^g
\theta(v;\Omega)\right|_{v=K^{\zeta}}$ first appeared in \cite{F2} in a different context.}

\begin{theorem}{\rm (cf. \cite {KK})} Let $\tau_\mu=\tau_\mu(C^\alpha, \omega)$ be given by formula (\ref{taukk}). Then\\
(i) $\tau_\mu$ does not depend on either $\zeta$ or the choice of the cuts
in the homology classes $a_i, b_i$;\\
(ii) $\tau_0$ is a nowhere vanishing holomorphic function on the moduli space $\Hc_{g}^0$
of generic holomorphic 1-differentials, whereas $\tau_\mu$ for a non-trivial $\mu$ is defined locally 
up to a root of unity and may depend on the choice of parameters $\zeta_k$;\\
(iii) $\tau_\mu$ is a solution of (\ref{tau}).
\end{theorem}

We want to describe how the tau function transforms under the action of $\C^*$ and $Sp(2g,\Z)$ 
on the moduli spaces $\Hc_{g}^\mu$ of holomorphic 1-differentials of degeneracy type $\mu$.
 
\begin{lemma}
The tau function $\tau_\mu$ on the space $\Hc_{g}^\mu$ has the property
\be
\tau_\mu(C^\alpha, \e \omega)=\e^{2\left(2g-2+r-\sum_{k=1}^r\f{1}{m_k+1}\right)}\; \tau_\mu(C^\alpha, \omega)
\la{taue}
\ee
for any $\e\in\C^*$. In other words, $\tau_\mu$ is a homogeneous function
on $\Hc_{g}^\mu$ of degree $2\left(2g-2+r-\sum_{k=1}^r \f{1}{m_k+1}\right)$.
\end{lemma}

\begin{proof} 
It is easy to see that the difference between 
the tau functions
$\tau_\mu(C^\alpha, \omega)$ and $\tau_\mu(C^\alpha, \e \omega)$ in (\ref{taukk}) comes from the different choice
of natural local parameters $\zeta$ on $C-\{x_1,\dots,x_r\}$ and $\zeta_k$ near $x_k\in C$. 
As above, we have $\zeta^\e=\e\zeta$ and $\zeta_k^\e=\e^{\frac{1}{m_k+1}}\zeta$.
Substituting these parameters $\zeta_k^\e$ into (\ref{taukk}), we get Eq. (\ref{taue}).
\end{proof}

\begin{corollary}
For the tau function $\tau_\mu$ on the space $\Hc_{g}^\mu$ we have the identity
\be
\sum_{i=1}^{2g+r-1} z_i\,\int_{s_i}\f{S_B-S_\omega }{\omega}= -\pi\sqrt{-1}\left(2g-2+r-\sum_{k=1}^r\f{1}{m_k+1}\right)\;.
\la{Euler1}
\ee
\end{corollary}

\begin{proof} 
The homogeneity property (\ref{taue}) implies that
$$\sum_{i=1}^{2g+r-1} z_i\f{\del}{\del z_i}\log \tau_\mu(C^\alpha, \omega) = 2\left(2g-2+r-\sum_{k=1}^r\f{1}{m_k+1}\right)\;.$$
This immediately yields (\ref{Euler1}) due to the definition (\ref{tau}) of the tau function.
\end{proof}

The behavior of the tau function under the change of Torelli marking of $C$ is described in the following lemma: 

\begin{lemma}
Let two canonical bases $\alpha=\{a_i,b_i\}_{i=1}^g$ and $\alpha'=\{a'_i,b'_i\}_{i=1}^g$ in $H_1(C)$
be related by $\alpha'=\sigma\alpha$, where 
\be
\sigma=
\left(\begin{array}{cc} D & C\\
B & A \end{array}\right)\in Sp(2g,\Z).
\la{symtrans}\ee
Suppose that the moduli space $\Hb_g^\mu$ parametrizes holomorphic differentials with at least one simple zero.
Then we have on $\Hb_g^\mu$
\be
\f{\tau_\mu(C^{\alpha'}, \omega)}{\tau_\mu(C^\alpha, \omega)}= 
{\rm det}(C\Omega+D)^{24}\;.
\la{tautaut}
\ee
where ${\Omega}$ is the period matrix of the Torelli marked Riemann sutface $C^\alpha$.
\end{lemma}

\begin{proof}
To establish this transformation property, we use the explicit formula (\ref{taukk}). 
According to Lemma 6 of \cite{KK}, when $\omega$ has at least one simple zero one can always choose 
the cut system on $C$ in such a way that $Z=Z'=0$ in (\ref{taukk}). The change of basis $\alpha'=\sigma\alpha$ 
results in the following transformation of the prime form $E(x,y)$:
\be 
E'(x,y)= E(x,y)e^{\sqrt{-1} \pi 
{\mathcal A}^x(y) (C\Omega +D)^{-1} C ({\mathcal A}^x(y))^t}
\la{transE} 
\ee 
(cf. \cite{F2}, Eq. (1.20)). 
For the expression 
$${\mathcal C}(x)=\frac{1}{W(x)}
\left.\left(\sum_{i=1}^g\omega_i(x)\frac{\partial}{\partial v_i}\right)^g
\theta(v;\Omega)\right|_{v=K^x}$$
it is shown in \cite{F2}, Eq. (1.23), that 
\be
{\mathcal C}'(x)=\sigma ({\rm det}(C{\Omega} +D))^{3/2} \, e^{\sqrt{-1}\pi K^x (C{\Omega}+D)^{-1}C (K^x)^t}\,  {\mathcal C}(x)\;,
\la{transC}
\ee
where $\sigma$ is a root of unity of degree 8, and $K^x$ is the vector of Riemann constants (\ref{rc}). 
Substituting these formulae into (\ref{taukk}), we obtain the statement of the lemma.
\end{proof}

Recall that there is one-to-one correspondence between $\C^*$-homogeneous holomorphic functions on $\Hc_g$ of
degree $n$ and holomorphic sections of the $n$-th power $L^n$ of the tautological line bundle $L\to\Hc_g^\mu/\C^*$. 
Since $\H_g^\mu=\Ht_g^\mu/\C^*=\Hc_g^\mu/Sp(2g,\Z)\times\C^*$, combining  Lemmas 2 and 3 we see that
the function $\tau_\mu(C^\alpha,\omega)$ on $\Hc_{g}^\mu$ descends to a non-vanishing holomorphic
section $\tau_\mu$ of the line bundle $\lambda^{24}\otimes L^{-2\left(2g-2+r-\sum_{k=1}^r \f{1}{m_k+1}\right)}$ on $\H_g^\mu=\Ht_g^\mu/\C^*$.\footnote{For a generic $\mu$ this may only be true for some power of $\tau_\mu$, cf. 
Theorem 1, (ii).} 
As a consequence we have
\begin{lemma}
In ${\rm Pic}(\H_g^\mu)\otimes\Q$ the following relation holds:
$$24\lambda-2\left(2g-2+r-\sum_{k=1}^r \f{1}{m_k+1}\right)\psi=0,$$
where $\psi=c_1(L)$.
\label{sect}\end{lemma}

\section{Divisor of the tau function}

Here we describe asymptotics of the tau function $\tau_0$ near boundary
components $D_{\deg}$, $D_0$ and $D_1,\dots,D_{[g/2]}$ of the space $\Hb_g$.
We start with the divisor $D_{\deg}$. In this case we may assume that the curve $C$ is fixed and consider a family $\omega_t$
of 1-differentials on $C$ such that two its simple zeroes, say, $x_{2g-3}(t)$ and $x_{2g-2}$(t), coalesce as $t\to 0$. We can take
a parametrization such that $t=(z_{4g-3}(t))^2=\left(\int_{x_{2g-2}(t)}^{x_{2g-3}(t)}\omega_t\right)^2$, where the integration path 
is chosen so that $t=0$ on $D_{\deg}$
(note that $t$ does not depend on a labeling of zeroes of $\omega$).
\begin{lemma}
The tau function has the following asymptotics near $D_{\deg}$:
\begin{equation}
\tau_0(C^{\a},\omega_{t}) =  t^{1/3 } \tau_1(C^{\alpha},\omega_0) (c+ o(1))
\label{taucaustic}
\end{equation}
for some constant $c\neq 0$.
\label{deg}\end{lemma}

\begin{proof}
Since the curve $C$ does not change, the bidifferential $B$ is independent of $t$, and we have
$$
\f{S_B- S_{\omega_t}}{\omega_t}\longrightarrow \f{S_B- S_{\omega_{0}}}{\omega_{0}}\quad {\rm as}\quad t\to 0
$$
or, equivalently, $\f{\del}{\del z_k(t)}\log\tau_0(C,\omega_t)\to\f{\del}{\del z_k(0)}\log\tau_1(C,\omega_0),\;k=1,\dots,4g-4.$
Therefore,
\be
\f{\tau_0(C^{\a},\omega_t)}{\tau_1(C^{\a},\omega_{0})}= c(t)(1+o(1))\quad {\rm as}\quad t\to 0, 
\la{hom1}
\ee
where $c(t)$ is independent of $z_1(t),\dots,z_{4g-4}(t)$.
To find $c(t)$ we use the homogeneity property (\ref{taue}) of the tau function.
Namely, according to (\ref{taue}), the homogeneity degree of the function $\tau_0(C^{\a},\omega_t)$ is $6g-6$ 
since here $r=2g-2$ and all $m_k=1$. 
By the same formula the degree of the function $\tau_1(C^{\a},\omega_{0})$ is  $(6g-6)-2/3$ since in this case 
$r=2g-3$, one of $m_k$ is equal to 2, and all other are equal to $1$.
The homogeneity degree of the local parameter $t$ is $2$, and, therefore, we have  $c(\e^2 t)/c(t)= \e^{2/3}$ for any $\e$. Thus, $c(t)= t^{1/3} c$ with some constant $c\neq 0$.
\end{proof}

The asymptotic of $\tau_0$ near the divisors $D_1,\dots, D_{[g/2]}$ can be computed in a similar way. Consider a
family $(C^\a_t,\omega_t),\;t\to 0,$ such that the limit curve $C^\a_0$ is a reducible curve with components $C_1^{\a_1}$ (of genus $j$) and
$C_2^{\a_2}$ (of genus $g-j$), where $\a=\a_1\cup\a_2$, and $\omega_t$ converges to $\omega_1$ on $C_1^{\a_1}$ and $\omega_2$ on $C_2^{\a_2}$,
where both $\omega_1$ and $\omega_2$ have only simple zeroes.
Under such a degeneration two simple zeroes of $\omega_t$, say,  $x_{2g-3}(t)$ and $x_{2g-2}$(t), tend to the node of $C_0^\a$. 
Again we may assume that $t=(z_{4g-3}(t))^2=\left(\int_{x_{2g-2}(t)}^{x_{2g-3}(t)}\omega_t\right)^2$, where the integration path gets
contracted to the node as $t\to 0$.

\begin{lemma}
At the limit $t\to 0$ the tau function $\tau_0$ has the following asymptotics near the boundary component $D_{j},\;j=1,\dots,[g/2]$:
\be
\tau_0(C^{\a}_t,\omega_t)=  t^3 \tau_0(C^{\a_1}_1, \omega_1)\,\tau_0(C^{\a_2}_2, \omega_2)(c+o(1)) 
\la{astaubound}
\ee
for some constant $c\neq 0$.
\label{red}\end{lemma}

\begin{proof}
First, we notice that away from an arbitrary neighborhood of the node we have
$$\f{S_{B_t} - S_{\omega_t}}{\omega_t}\to\left\{\begin{array}{cc}\f{S_{B_{1}}- S_{\omega_1}}{\omega_1}&{\rm on}\; C^{\a_1}_1\medskip\\  
\f{S_{B_{2}}- S_{\omega_2}}{\omega_2}&{\rm on}\; C_2^{\a_2}\end{array}\right.$$ 
(here $S_{B_{i}}$ is the Bergman projective connections on $C_{i}^{\a_{i}},\;i=1,2$).
Therefore, as in the previous lemma,
\be
\f{\tau_0(C^{\a}_t,\omega_t)}{\tau_0(C^{\a_1}_1,\omega_1)\tau_0(C^{\a_2}_2,\omega_2)}=c(t)(1+o(1)). 
\la{hom3}
\ee
Once again, to explicitly find $c(t)$ we use the homogeneity property (\ref{taue}) of the tau function.  
As we already know, the homogeneity degree of $\tau_0(C^{\a}_t,\omega_t)$ is $6g-6$,
whereas the degres of the functions  $\tau_0(C^{\a_1}_1,\omega_1)$ and $\tau_0(C^{\a_2}_2, \omega_2)$ equal $6j-6$ and $6(g-j)-6$ respectively.
Since the degree of the local parameter $t$ is $2$, this yields  $c(\e^2 t)/c(t)= \e^{6}$, so that  $c(t)=t^3\,c$ 
with some constant $c\neq 0$.
\end{proof}

Let us now describe the behaviour of the tau function at the boundary component $D_0$. This case is different from 
the other two considered above, and here we will follow the approach of \cite{F}. Take a family $(C^\a_t,\omega_t)$ 
such that $C_t\to C_0$ as $t\to 0$, where $C_0$ is an
irreducible curve with one node that we realize as a smooth genus $g-1$ curve $C'_0$ with two points $p$ and $q$ identified.
We assume that the cycle $a_g\in H_1(C_t)$ vanishes under the degeneration, so that $\{a_i,b_i\}_{i=1}^{g-1}$ is a 
canonical basis in $H_1(C'_0)$ (from now on we fix the homology bases and omit the superscript $\a$ that displays the dependence on the Torelli marking). 
Actually, we can take $t=e^{2\pi\sqrt{-1}B_g(t)/A_g(t)}$, where
$A_g(t)=\int_{a_g}\omega_t$ and $B_g(t)=\int_{b_g}\omega_t$. 
In particular, this means that $t\rightarrow 0$ as ${\rm Im}\,(B_g/A_g)\rightarrow \infty$.
Moreover, we can assume that the $a$-period $A_g(t)=A_g\neq 0$ is independent of $t$.
Then, in terms of the normalized
1-differentials $\omega_t=\sum_{i=1}^g A_i(t)\omega_i^t$, where $\omega_i^t$ converges to $\omega_i^0$, the $i$-th normalized
differential on $C'_0$, and $\omega_g^t$ tends to the meromorphic differential $\omega_{p,q}$ on $C'_0$ with simple poles at $p$ and $q$ with residues +1 and --1, and zero $a$-periods. Therefore, we have $\omega_t\to\omega_0=\sum_{i=1}^{g-1}A_i(0)\omega_i^0+(A_g/2\pi \sqrt{-1})\omega_{p,q}.$

\begin{lemma}
The tau function $\tau_0$ has the following asymptotics near the boundary component $D_0$: 
\be
\tau_0(C^{\a}_t,\omega_t)=  t^2 (c+o(1)),\quad\quad t\rightarrow 0,
\la{asimpD0}
\ee
where $c\neq0$ is a constant independent of $t$.
\label{irr}\end{lemma}

\begin{proof}
According to the definition of the tau function,
$$\f{\p}{\p B_g}\log \tau (C^{\a},\omega)=-\f{2}{\pi\sqrt{-1}}\int_{a_g}\f{S_B-S_\omega}{\omega}\;.$$
In the limit $t\rightarrow 0$, or, equivalently, ${\rm Im}\,(B_g/A_g)\to \infty$ we have
$$-\f{2}{\pi\sqrt{-1}}\int_{a_g}\f{S_B-S_\omega}{\omega}\longrightarrow - 4\,{\rm Res}_{\,p}\f{S_B^0-S_{\omega_0}}{\omega_0}\;.$$
From the definition of $S_{\omega_0}$ we immediately get
$${\rm Res}_{\,p}\f{S_B^0-S_{\omega_0}}{\omega_0}=\f{4\pi\sqrt{-1}}{A_g}\;.$$
Thus,
$$\f{\p}{\p B_g}\log \tau (C^{\a}_t,\omega_t)\longrightarrow \f{4\pi\sqrt{-1}}{A_g}\;$$
as $t\to 0$, which implies the asymptotics (\ref{asimpD0}).
\end{proof}

\begin{remark}
{\rm For $g=1$ the tau function is related to the Dedekind eta function
$\eta(\sigma)=e^{\pi i\sigma/12}\prod_{n=1}^\infty (1-e^{2\pi i n\sigma})$ by
$\tau(A,B)=\eta^{48}(B/A)$
(cf. \cite{KK}). The asymptotics (\ref{asimpD0}) obviously agrees with the asymptotics of the function $\eta$
as ${\rm Im}\,(B/A)\to\infty$.}
\end{remark}

Now we can prove the main result of this paper.

\begin{theorem}
In the rational Picard group ${\rm Pic}(\Hb_g)\otimes\Q$ of the space $\Hb_g=\P(\Eb_g)$ the following relation holds:
\be
\lambda=\f{g-1}{4}\psi+\f{1}{24}\d_\deg+\f{1}{12}\d_0+\f{1}{8}\sum_{j=1}^{[g/2]}\d_j\;.\label{mf}
\ee
Here $\lambda$ is the pullback of the Hodge class on $\Mb_g$ via the projection $\P(\Eb_g)\to \Mb_g$, $\psi$ is the tautological class 
on $\P(\Eb_g)$, $\d_\deg$ is the class of the divisor of degenerate 1-differentials, and $\d_j,\;j=0,\dots,[g/2],$ are the pullbacks of the classes of boundary divisors on $\Mb_g$.
\label{mt}\end{theorem}

\begin{proof}
Consider the divisor of the tau function $\tau_0$ on the space $\Hb_g$. From the above lemmas we know that it is supported on the
boundary $\Hb_g-\H_g$, so we have to compute the multiplicites of its components $D_\deg,D_0,\dots,D_{[g/2]}$. Let us start with
$D_\deg$. Choose a coordinate $\zeta$ transversal to $D_\deg$ such that $D_\deg$ is locally given by $\zeta=0$, and compute the
degree of the map $\zeta\mapsto t$, where $t=\left(\int_{x_{2g-2}(t)}^{x_{2g-3}(t)}\omega_t\right)^2$ is the parameter introduced in Lemma \ref{deg}. This is a local universal problem
well known in singularity theory. The double zero of a differential can be resolved in exactly three non-equivalent ways
(in accordance with the three decompositions of a permutation cycle of length 3 into the product of two transpositions).
Therefore, the map $\zeta\mapsto t$ is 3 to 1 for small $\zeta\neq 0$,
that is, $t=c\zeta^3+O(\zeta^4)$ as $\zeta\to 0$. This means that $\zeta=O(t^{1/3})$ as $t\to 0$, and from Lemma \ref{deg}, Eq.~(\ref{taucaustic}), it follows that the multiplicity of $D_\deg$ in the divisor of $\tau_0$ is 1.

The multiplicities of the divisors $D_0,\dots, D_{[g/2]}$ can be computed using similar considerations. 
Let us start with $D_0$. The parameter $t=e^{2\pi\sqrt{-1}B_g/A_g}$ introduced before Lemma \ref{irr} is a local
transversal coordinate to $D_0$, such that $D_0$ is locally given by the equation $t=0$. By Lemma \ref{irr},
Eq.~(\ref{asimpD0}), the multiplicity of $D_0$ in the divisor of $\tau_0$ is 2. For the divisors $D_2,\dots,D_{[g/2]}$
the parameter $t=\left(\int_{x_{2g-2}(t)}^{x_{2g-3}(t)}\omega_t\right)^2$ is a local transversal coordinate as well.
Therefore, by Lemma \ref{red}, Eq.~(\ref{astaubound}), each of these divisors enter the divisor of $\tau_0$ with multiplicity 3. The divisor $D_1$ requires a little more attention. Let $\zeta$ be a transversal coordinate
to $D_1$ such that $D_1$ is locally given by the equation $\zeta=0$. Then the natural map $\zeta\mapsto t$ is 2 to 1 
for small $\zeta\neq 0$. This happens because two pairs of diferentials $(\omega_1,\omega_2)$ and
$(-\omega_1,\omega_2)$ on a reducible curve $C_1\cup C_2$ with $g(C_1)=1$ and $g(C_2)=g-1$ represent the same point
in $D_1$, but are limits of two different families of differentials on smooth curves. Thus, the multiplicity of $D_1$
in the divisor of $\tau_0$ is 3/2, and we take care about the extra factor 1/2 by putting $\delta_1=1/2[D_1]$ (i.e.
considering $\Hb_g^0$ as a stack; note that $\delta_j=[D_j]$ for $j\neq1$).

Putting these computations together and recalling Lemma \ref{sect}, we get
$$24\lambda-(6g-6)\psi=\delta_\deg+2\delta_0+3\sum_{j=1}^{[g/2]}\delta_j,$$
which proves the theorem.
\end{proof}

\begin{remark} {\rm For a generic degeneracy type $\mu$ we can claim a somewhat less precise statement 
(the reason for that is the lack of detailed information about the irreducible components of the boundary for $\mu\neq 0$).
Namely, on the connected  components of the space $\Hb^\mu_g$ we have
\be
\lambda=\left(\f{g-1}{6}+\f{r}{12}-\f{1}{12}\sum_{i=1}^r\f{1}{m_i+1}\right)\psi + \delta, \label{kappa}
\ee
where $\mu=(m_1-1,\dots,m_r-1)$, and $\delta$ is an effective divisor supported on the boundary.
In this case the tau function $\tau_\mu$ is a holomorphic section of the line bundle 
$\lambda^{24}\otimes L^{-2\left(2g-2+r-\sum_{i=1}^r\f{1}{m_i+1}\right)}$ that is well-defined on each connected component
of $\Hb_g^\mu$ (this may in fact be true only for some integral
power of $\tau_\mu$, cf. Footnote 4 before Lemma \ref{sect}).}
\end{remark}

\section{Sums of the Lyapunov exponents}

Let us briefly review the Kontsevich-Zorich formula for the sum of the Lyapunov exponents of the $diag(e^t,e^{-t})$-action 
on the connected components of the space $\H_g^\mu$ (details can be found in the original article \cite{KZ}). The tangent space 
to the total space of the Hodge bundle $\E_g$ at a point $(C^\a,\omega)$ can be naturally identified with the relative homology group 
$H_1(C,\{x_1,\dots,x_r\},\C)$ by means of the period map associated with $\omega$ (cf. Section 3). The space $\E_g$ enjoys an invariant 
action of the group $GL_+(2,\R)$ that defines a real 4-dimensional foliation on $\E_g$. The foliation is $\C^*$-invariant and descends to
a 2-dimensional oriented foliation on $\P(\E_g)$ that preserves the stratification of $\P(\E_g)$ by the spaces $\H_g^\mu$. 
On each connected component $\M_g^\mu$ of the space $\H_g^\mu$ this 2-dimensional foliation is described by a closed form $\beta$ of real 
codimension 2. 
The Main Theorem of \cite{KZ} claims that for the sum $L_\mu=\lambda_1+\dots+\lambda_g$ of the Lyapunov exponents one has the following formula:
\be
L_\mu=\f{\int_{\M_g^\mu}\beta\wedge\lambda}{\int_{\M_g^\mu}\beta\wedge\psi}\;;
\label{kz}\ee 
here the Hodge class $\lambda$ and the tautological class $\psi$ are understood as elements of $H^2(\M_g^\mu,\Q)$.
It is conjectured in \cite{KZ} that $L_\mu$ is always rational.

The recent paper \cite{C} describes a construction of a ``Poincare dual'' to the form $\beta$ by means of Teichm\"uller curves.
More precisely, let $(T,t_0)$ be a once pointed elliptic curve. Consider the finite set of equivalence classes of branched covers 
of $T$ ramified only over $t_0$, of ramification type $(m_1+1,\dots,m_r+1,1,\dots,1)$ with $m$ entries equal to 1. For each such cover 
its degree $d$ and genus $g$ are related by the formulae $d=m_1+\dots+m_r+r+m$ and $d=2g-2+m+r$ (the former is the degree formula, 
and the latter is the Riemann-Hurwitz formula). The pullback of the differential $dz$ on $T$ to the cover has exactly
$r$ zeroes of multiplicities $m_1,\dots,m_r$. By changing the complex structure on $T$ one gets a complex curve 
(one dimensional Hurwitz space) that is a branched cover of the moduli space $\M_{1,1}$. Denote by $T_{d,\mu}$ its compactification 
in the sense of admissible covers. The connected components of $T_{d,\mu}$ naturally embed into $\M_g^\mu$ as complex 1-dimensional 
$SL(2,\R)$-invariant subvarieties called the {\em Teichm\"uller curves}. 
According to \cite{C}, \cite{EKZ} one has
\be
L_\mu=\lim_{d\to\infty}\f{T_{d,\mu}\cdot\lambda}{T_{d,\mu}\cdot\psi}\;.
\label{chen}\ee
Another relevant fact from \cite{EKZ} is that
$L_\mu=\kappa_\mu+c_\mu$, where 
$$\kappa_\mu=\f{g-1}{6}+\f{r}{12}-\f{1}{12}\sum_{i=1}^r\f{1}{m_i+1}$$
and $c_\mu$ is the Siegel-Veech area constant. Note that  the denominators in both (\ref{kz}) and (\ref{chen})
are closely related to each other and are relatively well understood -- the integral 
$\int_{\M_g^\mu}\beta\wedge\psi$ is essentially the volume of $\M_g^\mu$ \cite{EO},
and the intersection number $T_{d,\mu}\cdot\psi$ is the degree of the branched cover $T_{g,\mu}\to\Mb_{1,1}$ \cite{C}.

Our formula (\ref{kappa}) applied to (\ref{chen}) immediately yields
\be
L_\mu=\f{g-1}{6}+\f{r}{12}-\f{1}{12}\sum_{i=1}^r\f{1}{m_i+1}+\lim_{d\to\infty}\f{T_{d,\mu}\cdot\delta}{T_{d,\mu}\cdot\psi}\;,
\ee
so we recover the coefficient $\kappa_\mu$ and get an interpretation of the Siegel-Veech constant as a boundary
term. For the moduli space of generic differentials (that is, $\mu=0$) this interpretation can be made more precise. 
First, we observe that $T_{d,0}\cdot\delta_\deg=0$ by the construction of the Teichm\"uller curve $T_{d,0}$
(in this case $\omega$ has exactly $2g-2$ simple zeroes on any covering curve over any point in $\Mb_{1,1}$).
Second, a simple homological consideration shows that $T_{d,0}\cdot\delta_j=0$ for $j=1,\dots,[g/2]$
(an admissible cover of the degenerate elliptic curve remains connected after removing one node, so it cannot
represent a point in $D_j$ with $j\neq 0$). As it now follows from (\ref{mf}),
\be
L_0=\f{g-1}{4}+\f{1}{12}\lim_{d\to\infty}\f{T_{d,0}\cdot\delta_0}{T_{d,0}\cdot\psi}\;.
\label{zero}\ee
The numbers $T_{d,0}\cdot\delta_0$ are computable in the form of combinatorial sums
due to the transparent geometric construction of the Teichm\"uller curves $T_{d,0}$. Thus, in the case $\mu=0$, 
Eq. (\ref{mf}) allows us to reproduce the combinatorial formula for the Siegel-Veech constant $c_\mu$ (cf.
\cite{C}, \cite{EKZ}),
as well as the slope formula from \cite{C}, Section 3.

\section*{Acknowledgements} Both authors acknowledge the hospitality of the Max-Planck-Institut f\"ur Mathematik in Bonn 
and thank A.~Kokotov for helpful suggestions. 
In addition, PZ would like to express his gratitude to the Center for Quantum Geometry 
of Moduli Spaces at the \r{A}rhus University for support, and to thank M.~Kazarian and D.~Orlov for useful discussions. 
We are grateful to the referee for carefully reading the manuscript and proposing numerous improvements.

\end{document}